\RequirePackage{ifpdf}
\ifpdf 
\documentclass[pdftex]{sigma}
\else
\documentclass{sigma}
\fi

\def\s{{\mathfrak{S}}}

\def\ld{\lambda}

\def\ds{\displaystyle}

\def\ld{\lambda}

\begin{document}

\numberwithin{equation}{section}

\allowdisplaybreaks

\renewcommand{\PaperNumber}{072}

\FirstPageHeading

\renewcommand{\thefootnote}{$\star$}

\ShortArticleName{Skew Divided Dif\/ference Operators and Schubert Polynomials}

\ArticleName{Skew Divided Dif\/ference Operators\\ and Schubert Polynomials\footnote{This paper is a
contribution to the Vadim Kuznetsov Memorial Issue `Integrable
Systems and Related Topics'. The full collection is available at
\href{http://www.emis.de/journals/SIGMA/kuznetsov.html}{http://www.emis.de/journals/SIGMA/kuznetsov.html}}}

\Author{Anatol N.~KIRILLOV}

\AuthorNameForHeading{A.N.~Kirillov}

\Address{Research Institute of Mathematical Sciences (RIMS),
Sakyo-ku, Kyoto 606-8502, Japan}
\Email{\href{mailto:kirillov@kurims.kyoto-u.ac.jp}{kirillov@kurims.kyoto-u.ac.jp}}
\URLaddress{\url{http://www.kurims.kyoto-u.ac.jp/~kirillov/}}

\ArticleDates{Received May 01, 2007; Published online May 31, 2007}

 \Abstract{We study an action of the skew divided dif\/ference operators on the
Schubert polynomials and give an explicit formula
for structural constants for the Schubert polynomials in
terms of certain weighted paths in the Bruhat order on the symmetric group.
We also prove that, under certain assumptions, the skew divided
dif\/ference operators transform the Schubert polynomials into polynomials
with positive integer coef\/f\/icients.}


\Classification{05E15; 05E05}

\begin{flushright}
\it Dedicated to the memory of Vadim Kuznetsov
\end{flushright}

\section{Introduction}\label{section0}

In this paper we study the skew divided dif\/ference operators with
applications to the ``Little\-wood--Richardson problem'' in the Schubert
calculus. By the Little\-wood--Richardson problem  in the Schubert
calculus we mean
the problem of f\/inding a combinatorial rule for computing what one calls the
structural constants for Schubert polynomials.  These are the  structural
constants~$c^w_{uv}$, $u,v,w\in S_n$, of the ring $P_n/I_n$, where $P_n$
is the polynomial ring ${\bf Z}[x_1,\dots,x_n]$ and $I_n$ is the ideal
of $P_n$ generated by the symmetric polynomials without constant terms,
with respect to its ${\bf Z}$-free basis consisting of the classes of
Schubert polynomials~$\s_w$, $w \in S_n$.  Namely,
the constants $c^w_{uv}$ are
def\/ined via the decomposition of the product of two Schubert polynomials~$\s_u$ and $\s_v$ modulo the ideal $I_n$:
\begin{gather}
\s_v\s_u\equiv\sum_{w\in S_n}c_{uv}^w\s_w~({\rm mod}~I_n). \label{0.1}
\end{gather}
Up to now such a rule is known in the case when $u$, $v$, $w$ are the
Grassmannian permutations (see, e.g., \cite[p.~13]{M1} and \cite[Chapter~I,
Section 9]{M2})~-- this is the famous Littlewood--Richardson rule for Schur
functions~-- and in some special cases, see e.g.~\cite{Ko,Koh}.

The skew divided dif\/ference operators were introduced by I.~Macdonald in
\cite{M1}. The simplest way to def\/ine the skew divided dif\/ference operators is
based on the Leibniz rule for the divided dif\/ference operators $\partial_w$,
$w\in S_n$, namely,
\begin{gather}
\partial_w(fg)=\sum_{w\succeq
v}\left(\partial_{w/v}f\right)\partial_vg. \label{0.2}
\end{gather}
The symbol $w\succeq v$ for $w,v\in S_n$,
here and after, means that $w$ dominates $v$ with
respect to the Bruhat order on the symmetric group $S_n$ (see, e.g.,
\cite[p.~6]{M1}). Formula \eqref{0.2} is reduced to the classical Leibnitz rule
in the case when $w=(i,i+1)$ is a simple transposition:
\[
\partial_{i}(fg)=\partial_{i}(f)g+s_{i}(f)\partial_{i}g.
\]
One of the main applications of the skew divided dif\/ference operators
is an elementary and transparent algebraic proof of
the Monk formula for Schubert polynomials (see \cite[equation~(4.15)]{M1}).

Our interest to the skew divided dif\/ference operators is based on their
connection with the structural constants for Schubert polynomials.
More precisely, if $w,v\in S_n$, and $w\succeq v$, then
\begin{gather}
\partial_{w/v}\left(\s_u\right)|_{x=0}=c^w_{uv}. \label{0.3}
\end{gather}
The polynomial
$\partial_{w/v}(\s_u)$ is a homogeneous polynomial in $x_1,\ldots ,x_n$
of degree $l(u)+l(v)-l(w)$ with integer coef\/f\/icients. We make a
conjecture that in fact
\begin{gather}
\partial_{w/v}(\s_u)\in{\bf N}[x_1,\ldots ,x_n], \label{0.4}
\end{gather}
i.e.\ the polynomial $\partial_{w/v}(\s_u)$ has nonnegative integer
coef\/f\/icients.
In the case $l(u)+l(v)=l(w)$, this conjecture follows from the
geometric interpretation of the structural constants $c^w_{uv}$ as the
intersection numbers for Schubert cycles. For general $u,v,w\in S_n$ the
conjecture is still open.

In Section~\ref{section7} we prove conjecture \eqref{0.4} in the
simplest   nontrivial case (see Theorem~\ref{theorem1})
when $w$ and $v$ are connected by an edge in the Bruhat order on the
symmetric group $S_n$. In other words, if $w=vt_{ij}$, where $t_{ij}$ is
the transposition that interchanges $i$ and $j$, and $l(w)=l(vt_{ij})+1$.
It is well-known \cite[p.~30]{M1} that in this case the skew divided
dif\/ference operator $\partial_{w/v}$ coincides with operator $\partial_{ij}$,
i.e.\ $\partial_{w/v}=\partial_{ij}$.
Our proof employs the generating function for Schubert polynomials
(``Schubert expression'' \cite{FS,FK1,FK2}) in the nilCoxeter algebra.

In Section~\ref{section8} we consider another
application of the skew divided dif\/ference operators, namely, we give an
explicit
(but still not combinatorial) formula for structural constants for Schubert
polynomials in terms of weighted paths in the Bruhat order with weights
taken from the nilCoxeter algebra (see Theorem~\ref{theorem2}).

It is well known that there are several equivalent ways to def\/ine the skew Schur functions, 
see e.g., \cite{M1,M2}. Apart from the present paper, a few dif\/ferent def\/initions 
of skew Schubert polynomials have been proposed in 
\cite{BS,LS1} and \cite{CYY}. These def\/initions produce, in general,
 dif\/ferent polynomials.

\section{Skew Schur functions}\label{section1}

In this Section we review the def\/inition and basic properties of the skew
Schur functions. For more details and proofs, see \cite[Chapter I, Section 5]{M2}.
The main goal of this Section is to arise a problem of constructing skew
Schubert polynomials with properties ``similar'' to the those for skew
Schur functions (see properties \eqref{1.2}--\eqref{1.5} below).

Let $X_n=(x_1,\ldots ,x_n)$ be a set of independent variables, and $\lambda$,
$\mu$ be partitions, $\mu\subset\ld$, $l(\ld )\le n$.

\begin{definition} \label{definition1} The skew Schur function $s_{\ld /\mu}(X_n)$
corresponding to the skew shape $\ld -\mu$ is def\/ined to be
\begin{gather}
s_{\ld /\mu}(X_n)=\det\left(h_{\ld_i-\mu_j-i+j}\right)_{1\le i,j\le n},
\label{1.1}
\end{gather}
where $h_k:=h_k(X_n)$ is the complete homogeneous symmetric function of
degree $k$ in the variables $X_n=(x_1,\ldots ,x_n)$.
\end{definition}

Below we list the basic properties of skew Schur functions:

a) Combinatorial formula:
\begin{gather}
s_{\ld /\mu}(X_n)=\sum_Tx^{w(T)}, \label{1.2}
\end{gather}
where summation is taken over all semistandard tableaux $T$ of the
shape $\ld -\mu$ with
entries not exceeding $n$; here $w(T)$ is the weight of the tableau $T$
(see, e.g., \cite[p.~5]{M2}), and
$x^{w(T)}:=x_1^{w_1}x_2^{w_2}\cdots x_n^{w_n}$.

b) Connection with structural constants for Schur functions:
\begin{gather}
s_{\ld /\mu}=\sum_{\nu ,~l(\nu )\le n}c_{\mu\nu}^{\ld}s_{\nu}, \label{1.3}
\end{gather}
where the coef\/f\/icients $c_{\mu\nu}^{\ld}$ (the structural constants, or
the Littlewood--Richardson numbers) are def\/ined through the decomposition
\begin{gather}
s_{\mu}s_{\nu}=\sum_{\ld}c_{\mu\nu}^{\ld}s_{\ld}. \label{1.4}
\end{gather}

c) Littlewood--Richardson rule:
\begin{gather}
s_{\ld /\mu}=\sum_{\nu ,l(\nu )\le n} |{\rm Tab}^0(\ld -\mu ,\nu
)|s_{\nu}, \label{1.5}
\end{gather}
where ${\rm Tab}^0(\ld -\mu ,\nu )$ is the set of all semistandard
tableaux $T$ of shape $\ld -\mu$ and weight $\nu$ such that the
 reading word w$(T)$ of the tableaux $T$ (see, e.g., \cite[Chapter~I, Section 9]{M2}) is a lattice word (ibid). Thus,
\begin{gather}
{\rm Mult}_{V_{\ld}}(V_{\nu}\otimes V_{\mu})=c_{\nu\mu}^{\ld}=
|{\rm Tab}^0(\ld -\mu ,\nu )|. \label{1.6}
\end{gather}

\section[Divided difference operators]{Divided dif\/ference operators}\label{section2}

\begin{definition} \label{definition2} Let $f$ be a function of $x$ and $y$ (and possibly
other variables), the divided dif\/ference operators $\partial_{xy}$ is
def\/ined to be
\begin{gather}
\partial_{xy}f={f(x,y)-f(y,x)\over x-y}.\label{2.1}
\end{gather}
\end{definition}

The operator $\partial_{xy}$ takes polynomials to polynomials and has
degree $-1$. On a product $fg$, $\partial_{xy}$ acts according to the
Leibniz rule
\begin{gather}
\partial_{xy}(fg)=(\partial_{xy}f)g+(s_{xy}f)(\partial_{xy}g), \label{2.2}
\end{gather}
where $s_{xy}$ interchanges $x$ and $y$.

It is easy to check the following properties of divided dif\/ference
operators $\partial_{xy}$:
\begin{gather*}
\mbox{a)}\quad \partial_{xy}s_{xy}=-\partial_{xy},\qquad
s_{xy}\partial_{xy}=\partial_{xy},\\
\mbox{b)} \quad \partial_{xy}^2=0,\\
\mbox{c)} \quad \partial_{xy}\partial_{yz}\partial_{xy}=
\partial_{yz}\partial_{xy}\partial_{yz},\\
\mbox{d)} \quad \partial_{xy}\partial_{yz}=\partial_{xz}\partial_{xy}+
\partial_{yz}\partial_{xz}.
\end{gather*}

The next step is to def\/ine a family of divided dif\/ference operators
$\partial_i$, $1\le i\le n-1$, which act on the ring of polynomials
in $n$ variables.

Let $x_1,x_2,\ldots ,x_n$ be independent variables, and let
\[
P_n={\bf Z}[x_1,\ldots ,x_n].
\]
For each $i$, $1\le i\le n-1$, let
\[
\partial_i=\partial_{x_i,x_{i+1}},
\]
be the divided dif\/ference operator corresponding to the simple transposition
$s_i=(i,i+1)$ which interchanges $x_i$ and $x_{i+1}$.

Each $\partial_i$ is a linear operator on $P_n$ of degree $-1$. The divided
dif\/ference operators $\partial_i$, $1\le i\le n-1$, satisfy the following
relations
\begin{alignat*}{3}
& i)\quad && \partial_i^2=0, \ \ {\rm if} \ \ 1\le i\le n-1,&\\
& ii) && \partial_i\partial_j=\partial_j\partial_i, \ \ {\rm if}\ \ 1\le i,j\le
n-1, \ \ {\rm and} \ \ |i-j|>1,\\
& iii)\quad  && \partial_i\partial_{i+1}\partial_i=
\partial_{i+1}\partial_i\partial_{i+1}, \ \ {\rm if} \ \ 1\le i\le n-2.&
\end{alignat*}

Let $w\in S_n$ be a permutation; then $w$ can be written as a product of
simple transpositions $s_i=(i,i+1)$, $1\le i\le n-1$, namely,
\[
w=s_{i_1}\cdots s_{i_p}.
\]
Such a representation (or the sequence $(i_1,\dots,i_p)$) is called a
reduced decomposition of $w,$  if $p=l(w)$, where $l(w)$ is the
length of $w$. For each $w\in S_n$, let $R(w)$ denote the set of all
reduced decompositions of $w$, i.e.\ the set of all sequences $(i_1,\ldots
,i_p)$ of length $p=l(w)$ such that $w=s_{i_1}\cdots s_{i_p}$.

For any sequence ${\bf a}=(a_1,\ldots ,a_p)$ of positive integers, let us
def\/ine $\partial_{\bf a}=\partial_{a_1}\cdots\partial_{a_p}$.

\begin{proposition}[\cite{M1}, Chapter II]  \label{proposition1}\null {} \ \

i) If a sequence ${\bf a}=(a_1,\ldots ,a_p)$ is not reduced, i.e.\ not a
reduced decomposition of any $w \in S_n$, then
$\partial_{\bf a}=0$.

ii) If ${\bf a,b}\in R(w)$ then $\partial_{\bf a}=\partial_{\bf b}$.
\end{proposition}

From Proposition~\ref{proposition1}, $ii)$ follows that one can def\/ine
$\partial_w=\partial_{\bf a}$ unambiguously, where ${\bf a}$ is any
reduced decomposition of  $w$.

\section{Schubert polynomials}\label{section3}

In this section we recall the def\/inition and basic properties of the
Schubert polynomials introduced by A.~Lascoux and M.-P.~Sch\"utzenberger.
Further details and proofs can be found in~\cite{M1}.

Let $\delta =\delta_n=(n-1,n-2,\ldots ,1,0)$, so that
\[
x^{\delta}=x^{n-1}_1x^{n-2}_2\cdots x_{n-1}.
\]

\begin{definition}[\cite{LS}] \label{definition3} For each
permutation $w\in S_n$ the Schubert polynomial $\s_w$ is def\/ined to be
\begin{gather}
\s_w=\partial_{w^{-1}w_0}(x^{\delta}), \label{3.1}
\end{gather}
where $w_0$ is the longest element of $S_n$.
\end{definition}

\begin{proposition}[\cite{LS,M1}] \label{proposition2} \ \

i) For each permutation $w\in S_n$, $\s_w$ is a polynomial in $x_1,\ldots
,x_{n-1}$ of degree $l(w)$ with positive integer coefficients.

ii) Let $v,w\in S_n$. Then
\[
\partial_v\s_w=\begin{cases} \s_{wv^{-1}}, & if \ \   l(wv^{-1})=l(w)-l(v),\cr
0, & otherwise.\end{cases}
\]

iii) The Schubert polynomials $\s_w$, $w\in S_n$,
form a ${\bf Z}$-linear  basis in the space ${\cal F}_n$, where
\[
{\cal F}_n=\left\{ f\in
P_n~|~f=\sum_{\alpha\subset\delta}c_{\alpha}x^{\alpha}\right\} .
\]

iv) The Schubert polynomials $\s_w$, $w\in S_n$, form an orthogonal basis
with respect to the pairing $\langle ~,~\rangle_0$:
\[
\langle\s_w,\s_u\rangle_0
=\begin{cases} 1, & if \ \ u=w_0w,\cr 0, & otherwise,\end{cases}
\]
where by definition $\langle f,g\rangle_0=\eta
(\partial_{w_0}(fg)):=\partial_{w_0}(fg)|_{x=0},$ and $\eta(h)=h \vert_{x_1=
\cdots=x_n=0}$ for any polynomial $h$ in the variables $x_1,\dots,x_n.$

v) (Stability) Let $m>n$ and let $i:S_n\hookrightarrow S_m$ be the natural
embedding. Then
\[
\s_w=\s_{i(w)}.
\]
\end{proposition}

\section[Skew divided difference operators]{Skew divided dif\/ference operators}\label{section4}

The skew divided dif\/ference operators $\partial_{w/v}$, $w,v\in S_n$, were
introduced by I.~Macdonald~\cite[Chapter II]{M1}.

Let $w,v\in S_n$, and $w\succeq v$ with respect to the Bruhat order
$\succeq$ on the symmetric group~$S_n$. In other words, if ${\bf a}=(a_1,\ldots ,a_p)$ is a reduced
decomposition of $w$ then there exists a~subsequence ${\bf b}\subset{\bf
a}$ such that ${\bf b}$ is a reduced decomposition of $v$ (for more
details, see, e.g., \cite[equation~(1.17)]{M1}.

\begin{definition}[\cite{M1}] \label{definition4} Let $v,w\in S_n$, and $w\succeq
v$ with respect to the Bruhat order, and ${\bf a}=(a_1,\ldots a_p)\in
R(w)$. The skew divided dif\/ference operator $\partial_{w/v}$ is def\/ined
to be
\begin{gather}
\partial_{w/v}=v^{-1}\ds\sum_{{\bf b}\subset{\bf a},~{\bf b}\in R(v)}
\phi ({\bf a},{\bf b}), \label{4.1}
\end{gather}
where
\begin{gather*}
\phi ({\bf a},{\bf b})=\prod_{i=1}^p\phi_i({\bf a},{\bf b}),
\phi_i({\bf a},{\bf b})=\begin{cases} s_{a_i}, & if \ \ a_i\in{\bf b}, \cr
\partial_{a_i}, & if  \ \ a_i\not\in{\bf b} .\end{cases}
\end{gather*}
\end{definition}

One can show (see, e.g., \cite[p.~29]{M1}) that Def\/inition~\ref{definition4} is independent
of the reduced decomposition ${\bf a}\in R(w)$.

Below we list the basic properties of the skew divided dif\/ference
operators $\partial_{w/v}$. For more details and proofs, see, e.g., \cite{M1}.
The statement iv) of Proposition~\ref{proposition3} below seems to be new.

\begin{proposition} \label{proposition3} \ \

i) Let $f,g\in P_n$, $w\in S_n$, then
\begin{gather}
\partial_w(fg)=\sum_{w\succeq v}(\partial_{w/v}f)\partial_vg. \label{4.2}
\end{gather}

More generally,

ii) Let $f,g\in P_n$, $u,w\in S_n$, and $w\succeq u$ with respect to the the
Bruhat order. Then
\begin{gather}
\partial_{w/u}(fg)=\sum_{w\succeq v\succeq u}u^{-1}v(\partial_{w/v}f)
\partial_{v/u}g \label{4.3}
\end{gather}
(generalized Leibnitz' rule).

iii) Let $w=vt$, where $l(w)=l(v)+1$, and $t=t_{ij}$ is the transposition
that interchanges $i$ and $j$ and fixes all other elements of $[1,n]$.
Then
\begin{gather}
\partial_{w/v}=\partial_{ij}, \label{4.4}
\end{gather}
where $\partial_{ij}:=\partial_{x_ix_j}$.

iv) Let $w_0$ be the longest element of $S_n$. Then
\begin{gather}
w_0v\partial_{w_0/v}=\partial_{w_0v}. \label{4.5}
\end{gather}

v) Let $u,v,w\in S_n$, $w\succeq u$, and $l(w)=l(u)+l(v)$. Then
\begin{gather}
\partial_{w/u}\s_v=c^w_{uv}, \label{4.6}
\end{gather}
where $c^w_{uv}$ are the structural constants for the Schubert
polynomials $\s_w$, $w\in S_n$; in other words,
\[
\s_u\s_v\equiv\sum_{w\in S_n}c^w_{uv}\s_w\quad ({\rm mod}~I_n),
\]
where $I_n$ is the ideal generated by the elementary symmetric functions
$e_1(x_1,\ldots ,x_n), \ldots, e_n(x_1,$ $\ldots ,x_n)$.
\end{proposition}

\begin{proof} We refer the reader to \cite[p.~30]{M1} for proofs of statements
i)--iii).

iv) To prove the identity \eqref{4.5}, we will use the formula
\eqref{4.2} and the following result due to I.~Macdonald \cite[equation~(5.7)]{M1}):
\begin{gather}
\partial_{w_0}(fg)=\sum_{w\in S_n}\epsilon (w)\partial_w(w_0f)
\partial_{ww_0}(g), \label{4.7}
\end{gather}
where for each permutation $w\in S_n$, $\epsilon (w)=(-1)^{l(w)}$ is the
sign (signature) of $w$.

Using the generalized Leibnitz formula \eqref{4.2}, we can write the LHS~\eqref{4.7}
as follows:
\begin{gather}
\partial_{w_0}(fg)=\sum_{w_0\succeq v}v(\partial_{w_0/v}f)
\partial_v(g). \label{4.8}
\end{gather}
Comparing the RHS of \eqref{4.7} and that of \eqref{4.8}, we see that
\[
v(\partial_{w_0/v}f)=\epsilon (vw_0)\partial_{vw_0}(w_0f).
\]

To f\/inish the proof of equality \eqref{4.5}, it remains to apply the following
formula \cite[equation~(2.12)]{M1}:
\[
\partial_{w_0ww_0}=\epsilon (w)w_0\partial_ww_0.
\]

v) We consider formula \eqref{4.6} as a starting point for applications of the skew
divided dif\/ference operators to the problem of f\/inding a combinatorial
formula for the structural constants $c^w_{uv}$ (``Littlewood--Richardson
problem'' for Schubert polynomials, see Section~\ref{section1}). Having in mind some
applications of
\eqref{4.6} (see Sections~\ref{section7} and \ref{section8}), we reproduce below the proof of \eqref{4.6}
given by I.~Macdonald~\cite[p.~112]{M1}.
It follows from Proposition~\ref{proposition2} $ii)$ and Proposition~\ref{proposition3} $i)$ that
\[
c^w_{uv}=\partial_w(\s_u\s_v)=\sum_{w\succeq v_1}v_1
\left(\partial_{w/v_1}\s_u\right)\partial_{v_1}\s_v.
\]
In the latter sum the only nonzero term appears when $v_1=v$. Hence,
\[
c^w_{uv}=\partial_w(\s_u\s_v)=v\partial_{w/v}(\s_u)=\partial_{w/v}(\s_u),
\]
since~$\deg\partial_{w/v}(\s_u)=0$.
\end{proof}

It is well-known (and follows, for example, from Proposition~\ref{proposition2}, $i)$ and
$ii)$) that for each $v,w\in S_n$
\[
\partial_v\left(\s_w\right)\in{\bf N}[x_1,\ldots ,x_{n-1}].
\]
More generally, we make the following conjecture.

\begin{conjecture} \label{conjecture1} For any $u,v,w\in S_n$,
\[
\partial_{w/u}{\s_v}\in{\bf N}[X_n],
\]
i.e. $\partial_{w,u}(\s_v)$ is a polynomial in $x_1,\ldots ,x_n$ with
nonnegative integer coefficients.
\end{conjecture}

\begin{example} Take $w=s_2s_1s_3s_2s_1\in S_4$, $v=s_2s_1\in S_4$, and
${\bf a}=(2,1,3,2,1)\in R(w)$. There are three possibilities to choose
${\bf b}$ such that ${\bf b}\subset{\bf a}$, ${\bf b}\in R(v)$, namely, ${\bf
b}=(2,1,\cdot ,\cdot ,\cdot )$,\break ${\bf b}=(2,\cdot ,\cdot ,\cdot ,1)$
and ${\bf b}=(\cdot ,\cdot ,\cdot ,2,1)$. Hence,
\begin{gather*}
\partial_{w/v}=s_1s_2s_2s_1\partial_3\partial_2\partial_1+
s_1s_2s_2\partial_1\partial_3\partial_2s_1+
s_1s_2\partial_2\partial_1\partial_3s_2s_1 \\
\phantom{\partial_{w/v}}{} = \partial_3\partial_2\partial_1-\partial_1\partial_3\partial_{13}-
\partial_{13}\partial_2\partial_{14}.
\end{gather*}
Using this expression for the divided dif\/ference operator
$\partial_{w/v}$, one can f\/ind
\begin{gather*}
\mbox{a)} \quad \partial_{w/v}(x^3_1x^2_2)=x_1^2+x_1x_4+x_4^2\equiv x_2x_3 \ \ ({\rm
mod}~I_4).
\end{gather*}
Thus,
\[
\partial_{w/v}(x_1^3x_2^2)\equiv\s_{23}-\s_{13}+\s_{21} \ \ ({\rm
mod}~I_4).
\]
Here $\s_{23}$ means $\s_{s_2s_3},$ not $\s_{(2,4)}.$  Similar remarks
apply to other similar symbols here and after.

We used the following formulae for Schubert polynomials:
\begin{gather*}
\s_{23}=x_1x_2+x_1x_3+x_2x_3,\qquad \s_{13}=x_1^2+x_1x_2+x_1x_3,\qquad \s_{12}=x_1^2.
\\
\mbox{b)}\quad \partial_{w/v}(x_1^3x_2^2x_3)\equiv x_2^2x_3 \ \ ({\rm mod}~I_4),
\end{gather*}
 and
\begin{gather*}
\partial_{w/v}(x_1^3x_2^2x_3)\equiv\s_{121}+\s_{232}-\s_{123}-\s_{213}
-\s_{312} \ \ ({\rm mod}~I_4).
\\
\mbox{c)} \quad \partial_{13}(x_1^3x_2x_3)=x_1x_2x_3(x_1+x_3)\equiv
-x_1x_2^2x_3 \ \ ({\rm mod}~I_4).
\end{gather*}

Let us note that $x_1^3x_2x_3=\s_{12321}$.
\end{example}

These examples show that in general
\begin{itemize}\itemsep=0pt
\item the ``intersection'' numbers $\langle\partial_{w/v}(\s_u),
\s_{\tau}\rangle_0$ may  have negative values;

\item  coef\/f\/icients $c_{\alpha}$ in the decomposition
$\partial_{w/v}(\s_u)\equiv
 \sum\limits_{\alpha\subset\delta_n}c_{\alpha}x^{\alpha} \ \ ({\rm mod}~I_n)
$
may take  negative values.
\end{itemize}

\section[Analog of skew divided differences in the Bracket algebra]{Analog of skew divided dif\/ferences in the Bracket algebra}\label{section5}

In this Section for each $v,w\in S_n$ we construct the element $[w/v]$
in the Bracket algebra
${\cal E}_n^0$ which is an analog of the skew divided dif\/ference operators
$\partial_{w/v}$. The Bracket algebra ${\cal E}_n^0$ was introduced
in~\cite{FK2}. By
def\/inition, the Bracket algebra ${\cal E}_n^0$ (of type $A_{n-1}$) is the
quadratic algebra (say, over ${\bf Z}$) with generators $[ij]$, $1\le
i<j\le n$, which satisfy the following relations ${}$
\begin{alignat*}{3}
& (i) && [ij]^2=0, \quad {\rm for} \ \ i<j;&\\
& (ii) && [ij][jk]=[jk][ik]+[ik][ij],\quad
[jk][ij]=[ik][jk]+[ij][ik],\quad {\rm for} \ \ i<j<k; &\\
& (iii)\quad && [ij][kl]=[kl][ij]\quad {\rm whenever} \ \
\{ i,j\}\cap\{ k,l\} =\varnothing, \ \ i<j \ \ {\rm and} \ \ k<l.&
\end{alignat*}

For further details, see \cite{FK2,K1}.

Note that $[ij] \rightarrow \partial_{ij}$, $1 \le i < j \le n,$ def\/ines a
representation of the algebra ${\cal E}_n^0$ in $P_n.$

Now, let $v,w\in S_n$, and $w\succeq v$ with respect to the Bruhat order
on $S_n$. Let ${\bf a}\in R(w)$ be a~reduced decomposition of  $w$. We
def\/ine the element $[w/v]$ in the Bracket algebra ${\cal E}_n^0$ to be
\[
[w/v]=v^{-1}\ds\sum_{{\bf b}\subset{\bf a}, \ {\bf b}\in R(v)}
\phi ({\bf a},b),
\]
where
\[
\phi({\bf a},{\bf b})=\prod_i\phi_i({\bf a},{\bf b}),\qquad {\rm and}
\qquad \phi_i=\begin{cases} s_{a_i}, & a_i\in{\bf b},\cr [a_i~~a_i+1], &
a_i\not\in{\bf b}.\end{cases}
\]

Note that the right-hand side of the def\/inition of $[w/v]$ can be interpreted
inside the crossed product of ${\cal E}_n^0$ by $S_n$ (which is also called
a skew group algebra in this case) with respect to the action of $S_n$ on
${\cal E}_n^0$ def\/ined by
\[
w \cdot [w/v] = [w(i)w(j)]\qquad (\mbox{which means} \ \ -[w(j)w(i)] \ \ \mbox{if} \ \ w(i) > w(j)),
\]
eventually giving an element of ${\cal E}_n^0.$

\begin{remark} Let $w,v\in S_n$, and $w\succeq v$. One can show that the
element $[w/v]\in{\cal E}_n^0$ is independent of the reduced decomposition
${\bf a}\in R(w)$.
\end{remark}

\begin{conjecture} \label{conjecture2} The element $[w/v]\in{\cal E}_n^0$ can be written
as a linear combination of monomials in the generators $[ij]$, $i < j$,
with nonnegative  integer coefficients.
\end{conjecture}

\begin{example}
Take $w=s_2s_1s_3s_2s_1\in S_4$, $v=s_2s_1\in S_4$. Then
\begin{gather*}
[w/v]=[34][23][12]-[12][34][13]-[13][23][14]\\
\phantom{[w/v]}{}
=[34][12][13]+[34][13][23]-[12][34][13]-[13][23][14]\\
\phantom{[w/v]}{}=[13][14][23]+[14][34][23]-[13][23][14]=[14][34][23].
\end{gather*}
\end{example}

\section{Skew Schubert polynomials}\label{section6}

\begin{definition} \label{definition5} Let $v,w\in S_n$, and $w\succeq v$ with respect
to the Bruhat order. The skew
Schubert polynomial $\s_{w/v}$ is def\/ined to be
\begin{gather}
\s_{w/v}=\partial_{v^{-1}w_0/w^{-1}w_0}(x^{\delta_n}). \label{6.1}
\end{gather}
\end{definition}

\begin{example}
a) Let $w=s_1s_2s_3s_1\in S_4$, and $v=s_1\in S_4$.
Then $v^{-1}w_0=s_2s_1s_3s_2s_1$, $w^{-1}w_0=s_2s_1$, and
\begin{gather*}
\s_{w/v}=\partial_{21321/21}(x_1^3x_2^2x_3)=
(x_1^2+x_1x_4+x_4^2)x_2\\
\phantom{\s_{w/v}}{}\equiv\s_{121}+\s_{232}-\s_{124}
-\s_{213}-\s_{132} \ \ ({\rm mod}~I_4).
\end{gather*}

b) Take $w=s_3s_2\in S_4$ and $v=s_3\in S_4$. Then
$v^{-1}w_0=s_1s_2s_1s_3s_2$, $w^{-1}w_0=s_1s_2s_3s_2$, and
$\partial_{v^{-1}w_0/w^{-1}w_0}=\partial_{13}$. Thus
\[
\s_{w/v}=x_1^2x_2x_3(x_2+x_3)\equiv -x_1^3x_2x_3 \ \ ({\rm mod}~I_4).
\]

It is clear that if $w,v\in S_n$, and $w\succeq v$, then $\s_{w/v}$ is a
homogeneous polynomial of degree $\begin{pmatrix} n\cr 2\end{pmatrix} -l(w)+l(v)$ with
integer coef\/f\/icients. It would be a corollary of Conjecture~\ref{conjecture1} that skew
the Schubert polynomials have in fact positive integer coef\/f\/icients.
\end{example}

\begin{proposition} \label{proposition4} \ \

i) Let $v\in S_n$, and $w_0\in S_n$ be the longest element.
Then
\begin{gather}
\s_{w_0/v}=\s_{v}. \label{6.2}
\end{gather}

ii) Let $w\in S_n$, then
$\s_{w/1}=w_0ww_0\s_{ww_0}.$
\end{proposition}

Proof of \eqref{6.2} follows from \eqref{4.5} and \eqref{6.1}.

\medskip

It is an interesting task to f\/ind the Monk formula for skew Schubert
polynomials, in other words, to describe the decomposition of the product
$(x_1+\cdots +x_r)\s_{w/v}$, $w,v\in S_n$, $1\le r\le n-1$, in terms of
Schubert polynomials.

\section[Proof of Conjecture \ref{conjecture1} for divided difference operators
$\partial_{ij}$]{Proof of Conjecture~\ref{conjecture1} for divided dif\/ference operators
$\boldsymbol{\partial_{ij}}$}\label{section7}

First of all we recall the
def\/inition of the nilCoxeter algebra $NC_n$ and the construction of the
Schubert expression $\s^{(n+1)}\in{\bf N}[x_1,\ldots ,x_n][NC_n],$ where
${\bf N}[x_1,\ldots ,x_n][NC_n]$ denotes the set of all non-negative integral
linear combinations of the elements of the form $x_1^{m_1} \cdots x_n^{m_n}
\otimes e_{w},$ $m_1,\dots,m_n \in {\bf N}$, $w \in S_n,$ in ${\bf Z}[x_1,\dots,x_n] \otimes_{\bf Z} NC_n.$  Similar remarks apply to similar notation below.

 The study of action of divided dif\/ference operators
$\partial_{ij}$, $1\le i<j\le n$,
on the Schubert expression $\s^{(n+1)}$ is the main step of
our proof of Conjecture~\ref{conjecture1} for the skew divided dif\/ference operators
corresponding to the edges in the Bruhat order on the symmetric group
$S_{n+1}$. In exposition we follow to \cite{FS,FK1,FK2}.

\begin{definition} \label{definition6} The nilCoxeter algebra $NC_n$ is the algebra
(say, over ${\bf Z}$) with generators $e_i$, $1\le i\le n$, which satisfy
the following relations
\begin{enumerate}\itemsep=0pt
\item[$(i)$] $e_i^2=0$, for $1\le i\le n$,

\item[$(ii)$] $e_ie_j=e_je_i$, for $1\le i,j\le n$, $|i-j|>1$,

\item[$(iii)$] $e_ie_je_i=e_je_ie_j$, for $1\le i,j\le n$, $|i-j|=1$.
\end{enumerate}
\end{definition}

For each $w\in S_{n+1}$ let us def\/ine $e_w\in NC_n$ to be
$e_w=e_{a_1}\cdots e_{a_p}$, where $(a_1,\ldots ,a_p)$ is any reduced
decomposition of  $w$. The elements $e_w$, $w\in S_{n+1}$, are
well-def\/ined and form a ${\bf Z}$-basis in the nilCoxeter algebra $NC_n$.

Now we are going to def\/ine the Schubert expression $\s^{(n+1)}$ which
is a noncommutative generating function for the Schubert polynomials. Namely,
\begin{gather*}
\s^{(n+1)}=\sum_{w\in S_{n+1}}\s_we_w\in{\bf N}[x_1,\ldots ,x_n][NC_n].
\end{gather*}
The basic property of the Schubert expression $\s^{(n+1)}$ is that it
admits the following factorization~\cite{FS}:
\begin{gather}
\s^{(n+1)}=A_1(x_1)\cdots A_n(x_n), \label{7.1}
\end{gather}
where $A_i(x)=\prod\limits_{j=n}^i(1+xe_j)=(1+xe_n)(1+xe_{n-1})\cdots
(1+xe_i)$.

Now we are ready to formulate and prove the main result of this Section, namely,
the following positivity theorem:

\begin{theorem} \label{theorem1} Let $1\le i<j\le n+1$, $w\in S_{n+1}$. Then
\[
\partial_{ij}\s_w\in{\bf N}[x_1,\ldots ,x_{n+1}].
\]
\end{theorem}

\begin{proof} Our starting point is the Lemma below which is a generalization
of the Statement~4.19 from Macdonald's book~\cite{M1}. Before to state the Lemma,
we need to introduce a few notation.

Let $X=(x_1,\dots,x_n)$ be the set of variables and
$\mu=(\mu_1,\dots,\mu_p)$ be a composition of size $n.$ We assume that $\mu_j
 \ne 0$, $1 \le j \le p,$ and put by def\/inition $\mu_0=0.$

Denote by $X_j$ the set of variables $(x_{\mu_1+\cdots+\mu_{j-1}+1},\dots,
x_{\mu_1+\cdots+\mu_{j}}),$  $1 \le j \le p.$

Let $w \in S^{(n)}$ be a permutation such that the code of $w$ has
length $ \le n.$ The Schubert polynomial $\s_{w}(X)$ can be uniquely expressed
in the form
\[ \s_{w}(X)=\sum d_{u_1,\cdots,u_p}^{w}~\prod_{j=1}^{p} \s_{u_j}(X_j),
\]
summed over permutations $u_1 \in S^{(\mu_1)},\dots,u_p \in S^{(\mu_p)},$
see \cite[Chapter~IV]{M1}.

\begin{lemma} \label{lemma0} The coefficients $d_{u_1,\dots,u_p}^{w}$ defined above, are
non-negative integers.
\end{lemma}

The proof of Lemma proceed by induction on $l(u_p),$ and follows very close
to that given in~\cite{M1}. We omit details.

It follows from the Lemma above that it is enough to prove Theorem~\ref{theorem1}
only for the transposition  $(i,j)=(1,n)$.
Thus, we are going to prove that $\partial_{1n}\s_w\in{\bf N}[x_1,\ldots
,x_n]$. For this goal, let us consider the Schubert expression
$\s^{(n+1)}=A_1(x_1)A_2(x_2)\cdots A_n(x_n)$, see \eqref{7.1}. We are
going to prove that
\[
\partial_{1n}\s^{(n+1)}=\sum_{w\in S_{n+1}}\alpha_w(x)e_w,
\]
where $\alpha_w(x)\in{\bf N}[x_1,\ldots ,x_n]$ for all $w\in S_{n+1}$.
Using the Leibniz rule \eqref{2.2}, we can write
\begin{gather*}
\partial_{1n}\s^{(n+1)}=\partial_{1n}(A_1(x_1)A_2(x_2)\cdots A_n(x_n))\\
\phantom{\partial_{1n}\s^{(n+1)}}{}=\partial_{1n}(A_1(x_1))A_2(x_2)\cdots A_n(x_n)+A_1(x_n)A_2(x_2)\cdots
A_{n-1}(x_{n-1})\partial_{1n}(A_n(x_n)).
\end{gather*}

First of all,
\begin{gather*}
\partial_{1n}A_n(x_n)={1+x_ne_n-1-x_1e_n\over x_1-x_n}
=-e_n.
\end{gather*} The next observation is
\begin{gather*}
\partial_{1n}A_1(x_1)={A_1(x_n)-1\over x_n}+f(x_1,x_n),
\end{gather*}
where $f(x_1,x_n)\in{\bf N}[x_1,x_n][NC_n]$. Indeed, if
$A_1(x)=\sum\limits_{k=0}^{n}c_kx^k$, where $c_k\in NC_n$, $c_0=1$, then
\begin{gather*}
\partial_{1n}A_1(x_1)=\sum_{k=1}^{n}c_k{x_1^k-x_n^k\over x_1-x_n}=
\sum_{k=1}^{n}c_kx_n^{k-1}+f(x_1,x_n),
\end{gather*}
and $f(x_1,x_n)\in{\bf N}[x_1,x_n][NC_n]$, as it was claimed. Hence,
\begin{gather*}
 x_n\partial_{1n}\s^{(n+1)}=(A_1(x_n)-1)A_2(x_2)\cdots A_{n-1}(x_{n-1})
(1+x_ne_n)\\
\phantom{x_n\partial_{1n}\s^{(n+1)}=}{}-A_1(x_n)A_2(x_2)\cdots A_{n-1}(x_{n-1})e_nx_n+
x_nF(x_1,\ldots ,x_n)\\
\phantom{x_n\partial_{1n}\s^{(n+1)}}{}=A_1(x_n)A_2(x_2)\cdots A_{n-1}(x_{n-1})-A_2(x_2)A_3(x_3)\cdots A_n(x_n)\\
\phantom{x_n\partial_{1n}\s^{(n+1)}=}{}
+x_nF(x_1,\ldots ,x_n),
\end{gather*}
where $F(x_1,\ldots ,x_n)\in{\bf N}[x_1,\ldots ,x_n][NC_n]$. Thus, it is
enough to prove that the dif\/ference
\[
A_1(x_n)A_2(x_2)\cdots A_{n-1}(x_{n-1})-A_2(x_2)A_3(x_3)\cdots A_n(x_n)
\]
belongs to the set ${\bf N}[x_1,\ldots ,x_n][NC_n]$. We will use the
following result (see \cite{FS,FK1}):
\[
A_i(x)A_i(y)=A_i(y)A_i(x),\qquad 1\le i\le n.
\]
Thus, using a simple observation that $A_i(x)=A_{i+1}(x)(1+xe_i)$, we
have
\begin{gather*}
A_1(x_n)A_2(x_2)\cdots A_{n-1}(x_{n-1})=A_2(x_n)(1+x_ne_1)A_2(x_2)\cdots
A_{n-1}(x_{n-1})\\
\quad{}=A_2(x_n)A_2(x_2)\cdots A_{n-1}(x_{n-1})+x_nA_2(x_n)e_1A_2(x_2)\cdots
A_{n-1}(x_{n-1})\\
\quad{}=A_2(x_2)A_2(x_n)A_3(x_3)\cdots
A_{n-1}(x_{n-1})+x_nA_2(x_n)e_1A_2(x_2)\cdots A_{n-1}(x_{n-1})\\
\quad{}=A_2(x_2)A_3(x_n)(1+x_ne_2)A_3(x_3)\cdots A_{n-1}(x_{n-1})+
x_nA_2(x_n)e_1A_2(x_2)\cdots A_{n-1}(x_{n-1})\\
\quad{}=\cdots =A_2(x_2)A_3(x_3)\cdots A_n(x_n)+
x_n\sum_{i=1}^{n-1}\prod_{j=2}^i
A_j(x_j)A_{i+1}(x_n)e_i\prod_{j=i+1}^{n-1}A_j(x_j).
\end{gather*}

Let us denote  the sum over $i$ in (5.3) by $G(x_1,\ldots ,x_n)$. It is
clear that
\[
G(x_1,\ldots x_n)\in{\bf N}[x_1,\ldots x_n][NC_n].
\]
Thus the dif\/ference
\[
A_1(x_n)A_2(x_2)\cdots A_{n-1}(x_{n-1})-A_2(x_2)A_3(x_3)\cdots A_n(x_n)
=G(x_1,\ldots ,x_n)
\]
also belongs to the set ${\bf N}[x_1,\ldots ,x_n][NC_n]$.
\end{proof}

\section[Generating function for the Schubert polynomials structural
constants $c_{uv}^w$]{Generating function for the Schubert polynomials\\ structural
constants $\boldsymbol{c_{uv}^w}$}\label{section8}

Let $w,v\in S_n$, $l(w)-l(v)\le 1$, and $w\succeq v$ with respect to
the Bruhat order. For $1\le i\le n$ and $1\le s\le n-1$ we def\/ine the
element $e_i^{(s)}(w/v)$ of the nilCoxeter algebra
$NC_n$ using the following rule
\[
e_i^{(s)}(w/v)=\begin{cases} 0, & \mbox{if} \  \ w=vt_{(a,b)}, \  \ \mbox{and simultaneously} \  \
a\ne s \ \ \mbox{and} \  \ b\ne s,\cr e_{n-i}, &\mbox{if} \ \ w=vt_{(s,b)}, \ \ \mbox{and} \ \ s<b, \cr
-e_{n-i}, & \mbox{if} \ \ w=vt_{(b,s)}, \ \ \mbox{and}\ \ b<s,\cr
1, & \mbox{if} \ \  w=v.\end{cases}
\]

\begin{theorem} \label{theorem2} Let $u,w\in S_n$. Then
\begin{gather}
\sum_{v\in S_n}c_{uv}^we_v=\sum_{\{ v_i^{(s)}\}^{n-1}_{s=1}}
\prod_{s=1}^{n-1}\prod_{i=1}^{n-s}e_i^{(s)}\left(
v_{i-1}^{(s)}/v_i^{(s)}\right) , \label{8.1}
\end{gather}
summed over all sequences ${\bf v} =(v_i^{(s)})$ of permutations
\begin{gather*} w= v_0^{(1)}\succeq v_1^{(1)} \succeq\cdots\succeq v_{n-1}^{(1)}=v_0^{(2)}
\succeq v_1^{(2)} \succeq\cdots\succeq v_{n-2}^{(2)}\\
\phantom{w}{} =v_0^{(3)}\cdots=v_0^{(n-2)} \succeq v_1^{(n-2)} \succeq v_{2}^{(n-2)}=
v_0^{(n-1)} \succeq v_1^{(n-1)}=u
\end{gather*}
with restrictions
\[
l(v_{i-1}^{(s)})-l(v_{i}^{(s)}) \le 1 \qquad \mbox{for all $i$ and $s$.}
\]

In the product in the RHS \eqref{8.1} the factors are multiplied
left-to-right, according to the increase of $s$.
\end{theorem}

\begin{proof} We start with rewriting of the LHS \eqref{8.1}, namely showing that
\[
\sum_{v\in S_n}
c_{uv}^we_v=\eta\big(\partial_{w/u}\s^{(n)}\big),
\]
 where $\s^{(n)}$
denotes  the Schubert expression. Indeed,
\[
\eta\big(\partial_{w/v}\s^{(n)}\big)=\sum_{v\in S_n}\eta
\left(\partial_{w/v}(\s_v)\right) e_v=\sum_{v\in S_n}c_{uv}^we_v.
\]
The next step is to compute $\eta\left(\partial_{w/v}\s\right)$ using the
following lemma, which is obtained by repetitive use of the  generalized
Leibniz rule \eqref{4.3}.

\begin{lemma} \label{lemma1} Let $w,u\in S_n$, and $f_1,\ldots ,f_N\in P_n$. Then
\begin{gather*}
u\partial_{w/u}(f_1\cdots f_N)=\sum_{w=v_0\succeq v_1\succeq\cdots
\succeq v_{N-1}\succeq v_N=u}~\prod_{i=1}^Nv_i\left(\partial_{v_{i-1}/v_i}
(f)\right) .
\end{gather*}
\end{lemma}

We apply Lemma~\ref{lemma1} to the Schubert expression
\begin{gather}
\s^{(n)}=A_1(x_1)\cdots A_{n-1}(x_{n-1})=\prod_{i=1}^{n-1}\prod_{k=n-1}^i
(1+x_ie_k). \label{8.2}
\end{gather}
On the rightmost side of \eqref{8.2}, the factors are multiplied left-to-right
according to the increase of $i$.
As a result, we obtain
\[
\eta\big(\partial_{w/u}\s^{(n)}\big)=\sum_{\{ v_0^{(s)}\succeq v_1^{(s)}
\succeq\cdots\succeq v_{n-s}^{(s)}\}_{s=1}^{n-1}}\prod_{s=1}^{n-1}
\prod_{i=1}^{n-s}\eta\big(\partial_{v_{i-1}^{(s)}/v_i^{(s)}}
(1+x_se_{n-i})\big) ,
\]
summed over all sequences of permutations $\{ v_0^{(s)}\succeq
v_1^{(s)}\succeq\cdots\succeq v_{n-s}^{(s)}\}_{s=1}^{n-1}$ such that
$v_0^{(1)}=w$, $v_0^{(s+1)}=v_{n-s}^{(s)}$, $1\le s\le n-2$,
$v_1^{(n-1)}=u$.

Note that we omitted the action of the symmetric group elements since we
apply $\eta.$

It is clear that we can assume $l(v_{i-1}^{(s)})-l(v_i^{(s)})\le 1$ for
all $i$, $s$, and under these conditions, we have
\begin{gather*}
\eta\big(\partial_{v_{i-1}^{(s)}/v_i^{(s)}}(1+x_se_{n-i})\big)=
e_i^{(s)}\big(v_{i-1}^{(s)}/v_i^{(s)}\big).\tag*{\qed}
\end{gather*}
  \renewcommand{\qed}{}
\end{proof}

\section{Open problems}\label{section9}

Below we formulate a few problems related to the content of this paper.

{\bf 1. Main problem.} Let $w,v\in S_n$ and $w\succeq v$ with respect to the
Bruhat order on the symmetric group $S_n$. Prove that polynomials
$\partial_{w/v}\left(\s_u\right)$ have nonnegative coef\/f\/icients for each
$u\in S_n$.

{\bf 2. The generalized ``Littlewood--Richardson problem'' for Schubert
polynomials.}
Let $u,v,w\in S_n$ and
\[
\partial_{w/v}\left(\s_u\right)=\sum_{\alpha}c_{\alpha}x^{\alpha}.
\]
 Find a combinatorial description of the coef\/f\/icients
$c_{\alpha}:=c_{\alpha}(u,v,w)$.

\begin{remark} If $l(w)=l(u)+l(v)$ and $w\succeq v$, then $\partial_{w/v}
\left(\s_u\right)=c^w_{uv}$, see \cite[p.~112]{M1}, or the present paper,
Proposition~\ref{proposition3}, $v)$.
\end{remark}

{\bf 3. Skew key polynomials.} Let $\alpha =(\alpha_1,\ldots ,\alpha_n)$ be
a composition, $\lambda (\alpha )$ be a unique partition in the orbit
$S_n\cdot\alpha$, and $w(\alpha )\in S_n$ be the shortest permutation
such that $w(\alpha )\cdot\alpha =\lambda (\alpha )$. Let $v\in S_n$
be such that $w(\alpha )\succeq v$ with respect to the Bruhat order. Using
in Def\/inition~\ref{definition4} the isobaric divided dif\/ference operators
$\pi_i:=\partial_ix_i$, $1\le i\le n-1$  (see, e.g., \cite[p.~28]{M1}) instead
of operators $\partial_i$ one can def\/ine for each pair $w\succeq v$
 the skew isobaric divided dif\/ference operator $\pi_{w/v}:P_n\to P_n,$
where $P_n={\bf Z}[x_1,\dots,x_n].$

We def\/ine the skew key polynomial $k_{\alpha /v}$ to be
\[
k_{\alpha /v}=\pi_{w(\alpha )/v}\big( x^{\ld(\alpha )}\big) ,
\]
where $x^{\beta}=x_1^{\beta_1}\cdots x_n^{\beta_n}$ for any composition
$\beta =(\beta_1,\beta_2,\ldots ,\beta_n)$. It is natural to ask whether
or not the skew key polynomials have nonnegative coef\/f\/icients?

{\bf 4.} Find a geometrical interpretation of the skew divided dif\/ference
operators, the polynomials $\partial_{w/v}\left(\s_u\right)$, and the
skew key polynomials.

{\bf 5.} Does there exist a stable analog of the skew Schubert polynomials?

\subsection*{Acknowledgements}

This note is based on the lectures ``Schubert polynomials''
delivered in the Spring 1995 at the University of Minneapolis and in the
Spring 1996 at the University of Tokyo. I would like to thank my colleagues
from these universities
for hospitality and support.

The f\/irst version of this paper has appeared as a preprint \href{http://arxiv.org/abs/q-alg/9712053}{q-alg/9712053}.

\pdfbookmark[1]{References}{ref}
\LastPageEnding


\begin{thebibliography}{99}

\footnotesize\itemsep=0pt

\bibitem{BS} Bergeron N., Sottile F., Skew Schubert functions and
the Pieri formula for f\/lag manifolds,  {\it Trans. Amer. Math. Soc.} {\bf 354}
 (2002), no.~2, 651--673, \href{http://arxiv.org/abs/alg-geom/9709034}{alg-geom/9709034}.
 
 
 \bibitem{CYY} Chen W., Yan G.-G., Yang A., The skew Schubert 
polynomials,  {\it European J. Combin.}  {\bf 25}  (2004),  1181--1196.


\bibitem{FK1} Fomin S., Kirillov A.N., The Yang--Baxter equation,
symmetric functions and Schubert polynomials, {\it Discrete Math.}
{\bf 53} (1996), 123--143.

\bibitem{FK2} Fomin S., Kirillov A.N., Quadratic algebras, Dunkl
elements, and Schubert calculus, in Advances in Geometry, Editors
J.-L.~Brylinski and R.~Brylinski,  {\it Progr. Math.} {\bf 172} (1999), 147--182.

\bibitem{FS} Fomin S., Stanley R., Schubert polynomials and
nilCoxeter algebra, {\it Adv. Math.}  {\bf 103} (1994),  196--207.

\bibitem{K1} Kirillov A.N., On some quadratic algebras: Jucys--Murphy and
Dunkl elements,  in  Calogero--Moser--Sutherland models (1997, Montr\`eal, QC),  
{\it CRM Ser. Math. Phys.},  Springer, New York, 2000, 231--248, \href{http://arxiv.org/abs/q-alg/9705003}{q-alg/9705003}.


\bibitem{Ko} Kogan M., RC-graphs and a generalized Littlewood--Richardson
rule,  {\it Int. Math. Res. Not.}  {\bf 15} (2001), 765--782, \href{http://arxiv.org/abs/math.CO/0010108}{math.CO/0010108}.

\bibitem{Koh} Kohnert A., Multiplication of a Schubert polynomial by a
 Schur polynomial,  {\it Ann. Comb.} {\bf 1} (1997), 367--375.

\bibitem{LS} Lascoux A., Sch\"utzenberger M.-P.,
Polyn\^omes de Schubert, {\it C. R. Acad. Sci. Paris Ser. I Math.}
{\bf 294} (1982), 447--450.

\bibitem{LS1} Lenart C., Sottile F., Skew Schubert functions and the 
Pieri formula for f\/lag manifolds, {\it Trans. Amer. Math. Soc.}  {\bf 354} 
(2002), 651--673.

\bibitem{M1} Macdonald I.G., Notes on Schubert polynomials, Publ.
 LACIM,  Univ. du Quebec \`a Montr\'eal, 1991.

\bibitem{M2} Macdonald I.G., Symmetric functions and Hall
polynomials, 2nd ed., Oxford Univ. Press, New York, London, 1995.

\end{thebibliography}
\end{document}